\documentstyle[12pt,twoside]{amsart}

\title{Nash problem for stable toric varieities}
\author{Peter Petrov}
\date{April 19, 2006}
\address{Department of Mathematics, University of Georgia, Athens, GA 30602, USA
\newline
e-mail : petrov@@math.uga.edu}

\newcommand{\cone}{\operatorname{Cone}}
\newcommand{\sing}{\operatorname{Sing}}
\newcommand{\spec}{\operatorname{Spec}}
\newcommand{\cont}{\operatorname{Cont^n}}
\newcommand{\val}{\operatorname{val_v}}
\newcommand{\ord}{\operatorname{ord_t}}
\newcommand{\Hom}{\operatorname{Hom}}
\newcommand{\orb}{\operatorname{orb}}

\newcommand{\id}{\operatorname{id}}

\newtheorem{thm}{Theorem}[section]
\newtheorem{lem}[thm]{Lemma}

\newtheorem{defn}[thm]{Definition}
\newtheorem{Cor}[thm]{Corolary}
\newtheorem{Prop}[thm]{Proposition}
\theoremstyle{definition}
\newtheorem{Ex}[thm]{Example}

\begin{document}

\begin{abstract}
The main result of this paper is the proof of Nash conjecture for stable toric varieties. We also introduce Nash problem for pairs and prove it in the case of pairs $(X, Y$) of toric variety $X$ and a proper closed $T$-invariant subset $Y$, containing $\sing(X)$. 
\noindent
\end{abstract}

\maketitle
\markboth{\hfill Peter Petrov\hfill}{\hfill NASH PROBLEM FOR STABLE TORIC VARIETIES\hfill}
\section{Introduction}

The space of arcs $X_\infty$ of an algebraic variety $X$ encodes important information about its geometry. If $\pi\colon X_\infty \rightarrow X$ is the canonical map, then $\pi^{-1}(\sing(X))$ is a union of irreducible components. An irreducible component $C$ is called good if it contains an arc $\alpha$, such that $\alpha(\eta)\notin(\sing(X))$, where $\eta$ is the generic point of $\spec~k[[t]]$. In 1968 J. Nash proved that there is an injective map from the set of good components of $\pi^{-1}(\sing(X))$ to the set of the essential divisors over $X$, i.e. the exceptional divisors that "appear" on each resolution of singularities of $X$ (for the precise definition and main properties see Part 2). In the same paper \cite{Nash} Nash conjectured that this map is always bijective. Ishii and Koll\'ar proved in \cite{IK} that this is true for toric varieties, but fails in general, giving a counterexample in dimension 4, with one good component but two essential divisors. It was suggested that it would be useful to know the classes of algebraic varieties for which Nash conjecture holds (this is called Nash problem). Even in dimensions 2 and 3 the problem is open in general, although some results have been proved for ADE singularities (\cite{Nash}, \cite{Plenat}), minimal singularities (\cite{Reguera}), surface sandwich singularities (\cite{Reguera2}, \cite{LR}) and some other classes(\cite{PP}, \cite{Ishii1}). Recently Shihoko Ishii obtained a new result [private message] about Nash problem for toric pairs $(X, a)$, $X$ being toric variety and $a$ a toric ideal in $k[X]$. Although it is related to ours, our result does not follow from hers, because the class
of resolutions accepted in the former and the latter cases are different.

In this paper we give positive answer to Nash problem in the case of a stable toric variety (or STV) $X$. Because any such variety is a union of toric varieties, glued along toric subvarieties (see Part 3 for precise definitions), it is convenient to formulate the problem in the context of what we call Nash problem for pairs $(X, Y)$, where $X$ is an algebraic variety and $Y \subset X$ a proper nonempty closed subset (see Part 2). After we do it, we prove Nash problem for pairs in the case of a toric variety $X$ and a $T$-invariant closed subset $Y$. As a consequence, we obtain the main result for STV's.

The paper is organized as follows. In Part 2 we give the basic definitions and recall some facts that will be needed about Nash problem. Also we propose Nash problem for pairs as a new framework. This is of general interest as well. In Part 3 we briefly introduce STVs  and list some of their basic properties that will be needed later. In Part 4 the case of a toric pair $(X, Y)$ is considered, giving in Part 5 a positive answer to the Nash problem in the case of equidimensional STVs.\\
The autor would like to thank to Valery Alexeev for numerous suggestions and constuctive criticism, and to Shihoko Ishii and Willem Veys for the helpful remarks and comments.

\section{Nash problem}

In this paper varieties are always defined over an algebraically closed field $k$ of arbitrary characteristic. The scheme $\spec~k[[t]]$ has generic point $\eta$ and closed point $0$. A resolution of singularities of $X$ is a proper birational morphism $f:~Y \rightarrow X$ such that $Y$ is nonsingular and $f$ is an isomorphism on $Y\setminus f^{-1}(\sing(X))$.

The largest part of this section is based on \cite{IK} and \cite{Ishii}, which could be consulted for more motivation and details. 

Let $X$ be a scheme of finite type over $k$.
\begin{defn}
An $n$-jet of $X$ is a morphism $\spec~k[t]/(t^{n+1}) \rightarrow X$ and an arc is a morphism $\spec~k[[t]] \rightarrow X$.
\end{defn}
 
The set of all $n$-jets $X_n$ has structure of a scheme of finite type over $k$. In fact it represents the functor $F: Sch_k \rightarrow Set$, $F_n(Y):=\Hom_k(Y\times_{k}\spec~k[t]/(t^{n+1}, X)$. The space of arcs $X_\infty:=\varprojlim_{n}{X_n}$ is the projective limit of these schemes, which exist in the category of $k$-schemes because the transition morphisms are affine. It is a reduced scheme over $k$, in general not of finite type. For $m \geq n$ there is a morphism $\pi_m^n: X_m \rightarrow X_n$ induced by the truncation homomorphism $k[t]/(t^m) \rightarrow k[t]/(t^n)$ and a morphism $\pi_m^0: X_m \rightarrow X_0=X$. Also there are canonical morphisms  $\pi: X_{\infty}\rightarrow X$, $\pi(\alpha):=\alpha(0)$, and for each $n$, $\phi_n : X_{\infty} \rightarrow X_n$, corresponding to degree $n$ truncation of power series.\\
Now let $X$ be a $k$-variety with singular locus $\sing(X)$ and $\pi^{-1}(\sing(X))=\bigcup_{i\in I}C_i$ be the decomposition of a closed subscheme in $X_{\infty}$ into irreducible components.
\begin{defn}
$C_i$ is called a good component of $\pi^{-1}(\sing(X))$ (or for $X$) if it contains an arc $\alpha$ such that $\alpha(\eta) \notin \sing(X)$. 
\end{defn}
The next result is well known (\cite{IK}).
\begin{thm}
If $X$ is a toric variety, then every component is good. 
\end{thm}
This also holds in the case of arbitrary variety over a field with $char(k)=0$. But if characteristic $p\neq 0$ there may exist components which are not good (\cite{IK}).
\begin{defn}
Let $f_i: Y_i \rightarrow X, i=1,2$ be proper birational morphisms with normal $Y_i$'s and $E \subset Y_1$ be an irreducible exceptional divisor of $f_1$. Then $f_2^{-1}\circ f_1: Y_1 \dashrightarrow Y_2$ is defined on a nonempty open subset $E_0$ of $E$ and the center of $E$ in $f_2$ is defined to be  the closure of its image. We say that $E$ appears in $f_2$ if its center in $f_2$ is a divisor.
\end{defn}
This defines an equivalence relation $(E_1, f_1) \sim (E_2, f_2)$ on the set of pairs every time when $E_2$ is the center of $E_1$ in $f_2$. Each equivalence class corresponds to a divisorial valuation on $k(X)$ and is called exceptional divisor over $X$.
\begin{defn}
A resolution of singularities is called divisorial if the exceptional set is of pure codimension 1.
\end{defn}
\begin{Ex} For a $\mathbb{Q}$-factorial variety $X$, every resolution of singularities is divisorial.
\end{Ex}
\begin{defn}
An essential divisor of $X$ is an exceptional divisor $E$ over $X$ such that its center in each resolution of singularities of $X$ is an irreducible component of the exceptional locus. $E$ is divisorially essential if its center in every divisorial resolution is a divisor. 
\end{defn}
If $C_j$ is a good component for $X$, $\alpha \in C_j$ is an arc such that $\alpha(\eta) \notin \sing(X)$ and $f: Y \rightarrow X$ is any resolution of singularities, then by the valuative criterion of properness $\alpha$ lifts to a unique arc $\alpha ' \in Y_\infty$, such that $f \circ \alpha '= \alpha$. This gives $\alpha '(0) \in E$ for some exceptional divisor $E \subset Y$. Taking $\alpha$ to correspond to the generic point of $C_j$ (see \cite[Thm.2.15]{IK}), the divisor will be essential, so one can define Nash map for $X$, $\mathcal{N}_{\mathrm X}$, from the set of good components into the set of essential divisors.
\begin{thm}[Nash]
$\mathcal{N}_{\mathrm X}$ is always injective. In particular, the set of good components is finite.
\end{thm}
In the same paper \cite{Nash} Nash asked the following question:

\vspace{\stretch{1}}

{\bf Nash problem.} Is $\mathcal{N}_{\mathrm X}$ always bijective?

\vspace{\stretch{1}}

In their paper \cite{IK} Ishii and Koll\'ar proved that this fails in general. They give a counterexample with an affine hypersurface over an algebraically closed field of characteristic different from $2$ and $3$, which has one good component over its singular locus but two essential divisors. Their construction requires the dimension of the variety to be at least $4$, so the problem remains open in dimensions $2$ and $3$, although for some classes of surfaces a positive answer has been obtained (see the Introduction). Also, in the same paper the toric case is analyzed with a positive answer to Nash problem. 
\begin{thm}[Ishii, Koll\'ar]
For $X$ affine toric variety Nash map $\mathcal{N}_{\mathrm X}$ is bijective.
\end{thm}
Later in Section 4 we will modify the idea of their proof to obtain the main claim of this paper. For that, we will need some basic facts about the contact loci and the embedded version of Nash problem. 

Let $X$ be an affine variety over $k$, and $I$ be an ideal in the coordinate ring $k[X]$. Each arc $\alpha \in X_{\infty}$ defines homomorphism $\alpha^*: k[X] \rightarrow k[[t]]$. 
\begin{defn}
The $n$-th contact locus of $I$, $\cont(I)$, is the set of all $\alpha \in X_\infty$ for which $\min\{\ord_{t}~\alpha^*(f): f\in I\}=n$.
\end{defn}
These sets are important both for the singularity theory and for motivic integration. Any such set is a cylinder set, i.e. of the form $\phi^{-1}(S)$ for some constructible $S \subset X_n$ and some $n$. But its irreducible components are not cylinders in general (they are if $X$ is a smooth variety) and one could ask how to describe them.
\begin{defn}
For any discrete valuation of $k(X)$ associated with a prime divisor on some normal variety $X'$ birational to $X$, a positive integer times this valuation is called a divisorial valuation. 
\end{defn}
In \cite{ELM} the autors asked the following question: which divisorial valuations correspond to the irreducible components of contact loci? They call this the embedded version of Nash problem. In \cite{Ishii} Ishii gives an answer in the case of an affine toric variety $X$ corresponding to a cone $\sigma$, and $T$-invariant ideal $I \in k[X]$. In this case, every point $v \in \sigma \cap N$ corresponds to an orbit $T_{\infty}(v)$ under the action of $T_{\infty}$ on $X_{\infty}$ (\cite{Ishii}, 4.1). Also, any such point $v$ defines a divisorial valuation  by $\val_{v}(f):= \min \{(v, u): x^u ~\text{has non zero coefficient in}~ f \}$. Moreover, there is a partial order relation on the lattice points in $\sigma$:
\begin{defn} 
$v \leq_{\sigma} v'$ iff $v' \in v+ \sigma$.
\end{defn}
The next result (\cite{Ishii}, 4.12) gives us an important relation between the orbits corresponding to lattice points, in terms of this order in the case of an arbitrary toric variety $X=X(\Delta)$ defined by a fan $\Delta$.
\begin{defn}
For any face $\tau$ of $\sigma \in \Delta$, define $X(\tau)$ to be the closure of $\orb(\tau)$ and  
$X_{\infty}(\tau)$ to be $X(\tau)_{\infty} \setminus \bigcup_{\gamma~\text{not a face of}~\tau} X(\gamma)_{\infty}$.
\end{defn}
We denote by $U_{\gamma} \subset X$ the invariant affine open set containing only $\orb(\gamma)$ among the closed orbits, where $\gamma \in \Delta$. If $\tau < \tau' < \gamma$ are cones in $\Delta$, let $p: N_{\mathbb{R}}/ \tau \mathbb{R} \rightarrow  N_{\mathbb{R}}/ \tau' \mathbb{R}$ be the projection and let $\gamma':= p(\gamma)$.
\begin{thm}
If $T_{\infty}(v)$ and $T_{\infty}(v')$ are orbits in $X_{\infty}(\tau)$ and $X_{\infty}(\tau')$ respectively, then the next conditions are equivalent:

(i)~$\overline{T_{\infty}(v)}\supset T_{\infty}(v')$,

(ii)~$\tau < \tau'$ and there is a cone $\gamma > \tau'$ in the fan of $X$ such that $T_{\infty}(v), T_{\infty}(v') \subset (U_{\gamma})_{\infty}$ and $p(v) \leq_{\gamma'}v'$.
\end{thm}
The following theorem (\cite{Ishii}, 5.11) gives the answer to the embedded Nash problem for toric varieties. It will play an important role in the proof of the main result of this paper. 
\begin{thm}
The irreducible components of $\cont(I)$ are in bijection with the minimal elements of the set $V(I, n):= \{v \in \sigma \cap N: \min_{x^u \in I}(v, u)=n \}$ w.r.t. the partial order defined above. 
\end{thm}
Later, we will see that it becomes useful, in the case of STVs, to work with Nash problem for pairs \footnote{Suggested by V. Alexeev}, which we give now.

Let $(X, Y)$ be a pair of an algebraic variety $X$ and a proper closed subset $Y \supset \sing(X)$. 
\begin{defn}
A proper birational morphism $f: X' \rightarrow X$ with $X'$ smooth, such that $f^{-1}(Y)$ is of pure codimension 1, will be called a $Y$-resolution of $X$. The class of any prime divisor on $X'$ with the center appearing in any $Y$-resolution of $X$ as divisor will be called a $Y$-essential divisor over $X$. A good component of $\pi^{-1}(Y)$ is one that has an arc $\alpha$ such that $\alpha(\eta) \notin Y$.
\end{defn}
In the absolute case there is a difference between essential and divisorially essential divisors. Our definition is the analogue to the latter case. In \cite{IK} it is shown that for toric varieties the two notions coincide. In the case of toric pairs (see below) we will obtain a similar result in Sec.4.

Let $C$ be a good component for $(X, Y)$, $\alpha \in C$ is an arc such that $\alpha(\eta) \notin Y$ and $f': X' \rightarrow X$ is any $Y$-resolution of singularities. By the valuative criterion of properness, $\alpha$ lifts to a unique arc $\alpha ' \in X'_\infty$, such that $f \circ \alpha '= \alpha$. This gives $\alpha '(0) \in E'$ for some prime $Y$-exceptional divisor $E' \in X'$. 
\begin{thm}
For any pair $(X, Y)$ and any $Y$-resolution $X' \rightarrow X$ there is a map $${\mathcal{N}_{\mathrm (X, Y)}}: \{\text{good components of}~ (X, Y)\} \rightarrow \{\text{essential divisors of}~ (X, Y) \},$$ and it is injective.
\end{thm}
\begin{pf}[Sketch of the proof] The proof is virtually the same as the proof of Thm.2.15, \cite{IK}. For any $Y$-resolution $f: X' \rightarrow X$ and for any good component   $C$ of $(X, Y)$, if $z$ is its generic point, then by the remark above there exist some component of $\pi_{X'}^{-1}(f^{-1}(Y))$ (i.e. some irreducible component of the pre-image  of an $Y$-exceptional component $E$ on $X'_{\infty}$), whose generic point is sent by $f_{\infty}$ to $z$. To show that it is an essential component as well, we take another $Y$-resolution $X'' \rightarrow X$. Then by Def.2.4 $E'$ appears in $f''$, so by Def.2.16 it is $Y$-essential. This defines the map ${\mathcal{N}_{\mathrm (X, Y)}}$. It is an injective map, because if $C'$ is another good component and $z'$ is its generic point, then their lifts in $X'$ will be the generic point of the same $Y$-essential divisor. Taking its image by applying  $f'_{\infty}$ gives a generic point of one good component, which contradicts the choice of $z, z'$.
\end{pf}
We again call the map $\mathcal{N}_{\mathrm (X, Y)}$ Nash map. Then the next question arises naturally in this new contest.

{\bf Nash problem for pairs.} For which pairs $(X, Y)$ is the map $\mathcal{N}_{\mathrm (X, Y)}$ bijective?

In Section 4  we will prove that for a toric variety $X$ and a $T$-invariant closed subset $Y \supset \sing(X)$ the answer to this problem is positive. Such pairs will be called toric pairs.

\section{Stable toric varieties}
For more details about the definitions and results appearing in this section, see (\cite{Alexeev}).
\begin{defn}
A connected algebraic variety $X$ over $k$ (not necessarily irreducible) with action by a torus $T$ on $X$ is called a stable toric variety (STV) if it satisfies the following conditions:\\
i) X is seminormal;\\
ii) there are only finitely many orbits, and for each $x \in X$ the stabilizer $T_x \subset T$ is a subtorus.
\end{defn}
Stable toric varieties are analogs of stable curves in the case of toric varieties. Here we briefly give an idea for their classification. By \cite{Alexeev} each affine STV $X$ defines a face-fitting complex of cones $\Sigma$ with a reference map to $\Lambda_{\mathbb{R}}$, where ${\mathbb{Z}}^n \cong \Lambda \subset {\mathbb{R}}^n$ is a lattice. This means that we have a connected topological space $|\Sigma| = \cup \sigma_i$ and a finite-to-one map $\rho: |\Sigma| \rightarrow \Lambda_{\mathbb{R}}$ which identifies each $\sigma_i$ with a lattice cone. Since $\Sigma$ is face-fitting, the minimal faces of all $\sigma_i$ are equal to the same linear subspace $F_{\min} \subset \Lambda_{\mathbb{R}}$. Then every $\sigma _i$ is a preimage of a strictly convex cone in $X_{\mathbb{R}}/F_{\min}$. Moreover, $X$ is a union of (ordinary) toric varieties $X_{\sigma_{i}}$ glued in the way the complex $\Sigma$ is glued from $\sigma_i$, $X_{\sigma_i} \cap X_{\sigma_j} = \cup_{\sigma_i \cap \sigma_j=\tau}X_{\tau}$. Also, any projective polarized STV is glued from affine STVs in $\acute{e}$tale topology.
\begin{Ex} Let us take the complex of cones consisting of two cones in the plane, corresponding to the first and the third quadrants, with their faces. Then the corresponding STV will be two affine planes joined at the origin. Next, consider the complex which consist of the first and the second quadrants with their faces. Then the STV corresponding to it will be two planes intersecting along a line. But if we take the first quadrant (with all its faces) and a ray from the origin in say, the third quadrant, to form a complex, the corresponding STV will be a plane and a line intersecting it transversally. The first two constructions give examples of equidimensional STVs.
\end{Ex}

\section{The case of toric pairs}
The goal of this section is to prove the analogue of Thm.2.8 in the case of pairs. We note first that the analog of Thm.2.3 holds for pairs, so that one has:
\begin{lem}
For toric pair $(X, Y)$ all components of $\pi^{-1}(Y)$ are good.
\end{lem}
\begin{pf}
First, if $f: X' \rightarrow X$ is an equivariant $Y$-resolution of $X$, then the induced morphism $f_{\infty}$ is surjective. This is so because for any arc $\alpha \in X_{\infty}$, $\alpha(\eta) \in \orb(\tau)$ for some $\tau$ in the fan, defining $X$. Because $f$ is equivariant, the preimage of $\orb(\tau)$ contains a product $\orb(\tau) \times T'$ for some torus $T'$ of dimension less than $\dim X$. Thus the restriction of $\alpha$ on $k((t))$ lifts to $Y$ which, by the valuative criterion of properness, lifts $\alpha$ itself.

Next, suppose that $C$ is not a good component for $(X, Y)$. For any equivariant $Y$-resolution $f$ as above, with $E_i$ the irreducible components of the exceptional locus, $\pi_{X'}^{-1}(E_i)$ are the irreducible components of $f_{\infty}^{-1}(\pi^{-1}(Y))$. By the argument above, there exists $i$ such that $\pi_{X'}^{-1}(E_i)$ will be mapped to $C$. But the preimage of $E_i$ contains an arc which sends $\eta$ outside $\sing(X)$, which contradicts the choice for $C$.
\end{pf}
\begin{thm}
Let $X$ be a toric variety over $k$ and let $Y \subset X$ be a $T$-invariant proper closed subset containing $\sing(X)$. Then Nash map $\mathcal{N}_{\mathrm (X, Y)}$ is bijective.
\end{thm}
Let us make some remarks. First, we want to relate somehow the sets of $Y$-irreducible components and $Y$-essential divisors using the combinatorial data which defines $X$. Obviously, the question is local, so without loss of generality we could take $X$ to be an affine toric variety, defined by a cone $\sigma \subset N_{\mathbb{R}}$. Next, denoting by $I_Y \subset k[X]$ the ideal corresponding to $Y$, we see that $\pi^{-1}(Y)= \cup_{n \geq 1}\cont(I_Y)$. By Thm. 2.14, if $O_1$, $O_2$ are orbits in $X_{\infty}$ corresponding to points $v_1, v_2 \in N \cap \sigma$, then $\overline{O_1} \supset O_2$ iff $v_1 \leq_{\sigma} v_2$. Since $Y$ is $T$-invariant, it corresponds to a finite union of faces  of $\sigma$. Let ${\tau_1,..., \tau_s} \subset \sigma$ be all faces such that $\orb(\tau_i) \subset Y$. Now denote $$W^{\geq 0}:= \{v \in N \cap \sigma: (v, u) \geq 0~\text{for all}~ x^u \in I_Y \}$$ and $$W^0:= \{v \in N \cap \sigma: \exists x^u \in I_Y ~\text{s.t.}~ (v, u)=0\}.$$ Then $W^{\geq 0} \setminus W^0$ will contain exactly the lattice points which correspond to the orbits in $X_{\infty}$ contained in $\pi^{-1}(Y)$. But this set is actually $\cup_{i} \tau_i^{\circ}$. Indeed, if $I_Y = (x^u)$ is principal, then the difference of sets is just $(\sigma \cap N) \setminus H_{u}$ where $H_u$ is the hyperplane in $N_{\mathbb{R}}$ defined by $u$. In general, if $I_Y= \{x^{u_1},...,x^{u_r}\}$, then $\cap_{i=1,...,r}(\sigma \setminus H_{u_i}) = \sigma \setminus (\cup H_{u_i})$, which is in fact $W^{\geq 0} \setminus W^0$.

Another remark we want to make is about subdivisions of the cone~$\sigma$. We are interested in regular subdivisions, that is, subdivisions into regular cones, corresponding to resolutions of singularities $f:X' \rightarrow X$, and such that $f^{-1}(Y)$ is a divisor. If $\mathcal{F}$ is the map of fans corresponding to $f$, then for each $\nu$ in the set of cones defining $Y$, ${\mathcal{F}}^{-1}(\nu)$ is a union of cones each having either a ray $\rho$ with $\orb(\rho) \subset Y$ or a ray $\rho' \subset {\mathcal{F}}^{-1}(\nu \setminus \cup_{\rho ~\text{ray in}~ \nu} \rho)$. In the latter case this means  that $\rho'$ is not a ray of $\nu$ but of the fan obtained by the subdivision $\mathcal{F}$.
\begin{defn}
Define $$W:=\{v \in N \cap (\cup_i \tau_i^{\circ}): v~ \text{is minimal w.r.t.}~ \leq_{\sigma}\} \subset N \cap \sigma,$$ where $\tau_i^{\circ}$ is the relative interior of the cone $\tau_i$ for each $i$.
\end{defn}
\begin{lem}
There exists an injection $F_1$ from $W$ to the set $\mathcal{C}$ of irreducible components of $\pi^{-1}(Y)$.
\end{lem}
\begin{pf} For each element of $w \in W$, there is an arc $\alpha$ such that $\alpha(0) \in Y$, $\alpha(\eta) \in T$ and $w \in N$ is defined by the ring homomorphism $\alpha^*: k[X] \rightarrow k[[t]]$, taking $\alpha^*(x^u):=t^{(w, u)}$ (see \cite[Props.3.10, 3.11]{IK}). Also, there is a face $\tau$ containing $w$, such that $\alpha(0) \in \orb(\tau) \subset Y$ because $(w, u) \geq 0$ for all $u \in {\tau}^*$, thus $\alpha^*$ extends to $U_{\tau}$. Then $\alpha$ defines a good component $C$ of $(X, Y)$, so we have a map $F_1: W \rightarrow {\mathcal{C}}$. The injectivity follows from the fact that there is a non-empty open subset $U \subset C$, such that for all $\gamma \in C$ the correspondong lattice point will be $w$.
\end{pf}
\begin{defn} 
Define toric $Y$-divisorially essential divisors  to be the divisors which appear (see Def.2.4) in each toric $Y$-resolution of $X$.
\end{defn}
The next lemma is a modification of Lem.3.15 in \cite{IK}, so we skip some minor details in the proof.
\begin{lem}
There is a map $$F_2: \{\text{toric $Y$-divisorially essential divisors}\} \rightarrow W,$$ defined by $F_2(D_v):= v$, which is injective.
\end{lem}

\begin{pf} In any $Y$-resolution of $X$ defined by a fan $\Sigma$, each exceptional divisor corresponds to a ray $\rho \in \Sigma$ which either subdivides some face $\tau$ of $\sigma$, among the faces corresponding to $Y$, or coincides with it (i.e. $\rho = \tau$). So it is defined either by some primitive vector $w \in \tau^{\circ}$ or by some primitive vector $w \in \rho^{\circ}$. The important part to prove is that if a primitive vector $w \in N \cap (\cup \tau^{\circ})$ is not minimal, then the corresponding divisor $D_w$ defined by it does not appear in some $Y$-resolution. For this we will construct a $Y$-regular subdivision $\Sigma$ of $\sigma$ (that is, one corresponding to a $Y$-resolution) in which the ray $\rho = {\mathbb R}^{\geq 0}.w$ does not appear. 

Take such non-minimal $w$. Then there are $n_1, n_2 \in N \cap \tau$ such that $w = n_1+n_2$. It is easy to see that either
\begin{enumerate}
\item $n_1, n_2 \in W$, or
\item $n_1 \in W$, and $n_2$ could be taken on a ray of $\sigma$.
\end{enumerate}
This is so because if, say, $n_2$ is not in $N \cap (\cup \tau^{\circ})$, then $n_2$ will be in a non-singular face $\gamma$, generated by primitive vectors $p_1,...p_s$. As $n_2$ is a linear combination of $p_i$'s we have, say, a non-zero coefficient before $p_1$ in $n_2 = \sum_{i=1}^s b_{i}p_i$. Let $\delta$ is the minimal face, containing $n_1$ and $\sum_{i=2}^{s} a_{i}p_i$. But $\delta$ is singular (because contains $n_1$), and $n_1+ \sum_{i=2}^{s} a_{i}p_i \in \delta^{\circ}$. So replacing $n_1$ by $n_1+ \sum_{i=2}^{s} a_{i}p_i$ and $n_2$ by $a_{1}p_1$, we will have the case (2).

Take the minimal regular subdivision of $\cone (n_1, n_2)$ (it exists for any 2-dimensional cone)and let $\cone (u, v)$ be the cone in this subdivision, containing $w$ in its interior. Because at least one of $u, v$ should be in $W$, let $u \in W$. In the case (1) one takes the star subdivision $\Sigma'$ of $\Sigma$ with center $u$ and then takes the star subdivision of $\Sigma''$ of $\Sigma'$ with center $v$. The last subdivision then could be completed to a regular $Y$-subdivision, the fan of which would not contain the ray $\rho$. 

In the case (2), if $\Sigma'$ is not simplicial, take a minimal-dimensional cone $\mu$ with a lattice vector in its interior and the corresponding star subdivision of $\mu$. Continue this way to obtain a simplicial subdivision $\Delta_2$ with exceptional set of pure codimension 1. If it is not regular, take a cone $\beta$ with maximal multiplicity. The volume of the polytope $P= \sum a_{j}q_j$, generated by the primitive vectors $q_j$ on its rays, is then bigger than 1. This polytope contains a non-zero lattice point $m$ not lying on any of its edges, so one can take the star subdivision of $\beta$ with center $m$. Its exceptional set is a divisor and the volume of the corresponding polytopes will decrease or remain the same. Repeating this procedure for each cone with maximal multiplicity (bigger than 1) we will obtain a regular subdivision $\Delta_3 \supset cone(u, v)$ with exceptional locus of pure codimension 1. If necessary, we will perform few more subdivisions of the faces by rays ouside of $cone(u, v)$ to obtain at the end a $Y$-regular subdivision. Any of this subdivisions did not change $\cone(u, v)$, so $\rho \notin \Delta_3$. Also, all regular cones in $\Sigma$ did not change, so this will define the needed $Y$-resolution of $X$.

The injectivity is obvious by the definition of $F_2$.\\
\end{pf}
Thus we defined maps:
$$
\begin{array}{ccc}
$W$
&
  \stackrel{F_1}{\longrightarrow} 
&
\left\{
\begin{array}{c}
\mbox{the components of}\\
\mbox{ $\pi^{-1}(Y)$}
\end{array}
\right\} \\

&& \hphantom{{\cal N}}\downarrow {{\cal N}_{(X, Y)}}\\
{ F_2}\uparrow\hphantom{{ F_2}}
&
 
&
\left\{
\begin{array}{c}
\mbox{Y-essential divisors}\\
\mbox{of $X$}
\end{array}
\right\}\\
 & &  \mbox{ $\cap$}\\
\left\{
\begin{array}{c}
\mbox{toric Y-divisorially}\\
\mbox{essential divisors of X}\\

\end{array}
\right\} 
&
\supset
&
\left\{
\begin{array}{c}
\mbox{ Y-divisorially essential}\\
\mbox{ divisors of $X$}
\end{array}
\right\}\\
\end{array}
$$\\
The maps above satisfy the following
\begin{lem}
$F_2 \circ {\mathcal N}_{(X, Y)} \circ F_1 = {\id}_W$.
\end{lem}

\begin{pf} The map ${\mathcal N} \circ F_1: W \rightarrow \lbrace Y- \text{essential divisors of}~ X \rbrace$ sends a point $w$ to $D_w$. The reason is that the generic point of $F_1(w)$ is an arc, which could be lifted by a toric divisorial resolution $f: X' \rightarrow X$ to an arc $\gamma \in X'_{\infty}$ by the same argument as in the remark preceding Thm.2.7. So $\gamma(0)$ is the generic point on ${\mathcal N} \circ F_1 (w)$. Hence, the corresponding exceptional divisor defined by a ray $\rho$ contains $\gamma(0)$ and thus satisfies $\rho = D_w$. Then we apply Lem.4.5.
\end{pf}
Now we obtain the key result in this paper, Thm.4.2. The claim follows immedeately by applying Lem.4.4, 4.6 and 4.7.
Here are some consequences, coming directly from the diagram above.
\begin{Cor}
Given a toric variety $X$ with proper closed $Y \supset \sing(X)$, one has:
\\i) the set of the $Y$-essential divisors coincide with the set of divisorially $Y$-essential divisors over $X$;
\\ii) the number of $Y$-essential components is finite;
\\iii) the number of components of $\pi^{-1}(Y)$ is finite.
\end{Cor}

\section{Nash problem for STVs}
Now let $X$ be a stable toric variety, affine or polarized projective. As mentioned above, for our purposes without loss of generality one could take it to be affine. I.e., it corresponds to a complex $\Sigma$ of rational plyhedral cones. The singular locus of $X = X_{\Sigma}$ is a union of two sets: the union of intersection loci $\cup_{i \neq j} \lbrace X_{i} \cap X_{j} \rbrace$ and the union of singular loci $\cup_{i} \sing(X_{i})$. Moreover, there is a normalization map $\nu : \coprod_i X_i \rightarrow X_{\Sigma}$ such that for each fixed $i$, $\nu^{-1}(\sing(X)) \cap X_i = \sing(X_i) \cup (\cup_{j \neq i}(X_i \cap X_j))$. For each $i$ this gives a closed subset $Y_i \subset X_i$. Also, each essential divisor over $X_{\Sigma}$ becomes a $Y$-essentail divisor for the pair $(X_{\Sigma}, Y)$, where $Y:= \cup_i~ \nu(Y_i)$. Then the answer of Nash problem for $X$ comes naturally from the answers of Nash problem for each pair $(X_i, Y_i)$. This will follow from the Prop.5.1 and Prop.5.2 below.
Let $\Omega$ be a disjoint union of the sets of irreducible components of $\pi_i^{-1}(Y_i)$, where $\pi_i: {X_i}_{\infty} \rightarrow X_i$ .
\begin{Prop}
There is a one-to-one correspondence between the set of irreducible components of 
$\pi^{-1}(\sing(X))$ and $\Omega$.
\end{Prop}
\begin{pf}
Let $\alpha$ be an arc in $\pi^{-1}(Y)$, corresponding to the generic point of a component $C$ of $\pi^{-1}(\sing(X))$. Let $\nu_{\infty}$ be the map of arc spaces corresponding to the normalization map $\nu$. Then $\nu_{\infty}^{-1}(\alpha)$ contains the generic points of some components of $\pi_i^{-1}(Y_i)$ for some $i$'s. By the description of affine stable toric varieties preceding Def.3.3, we see that $i$ is unique. Taking the restriction of $\pi$ on the pre-image of $X_i$, we obtain that the component $C$ above is unique as well. This defines an injective map between the sets above. Conversely, take any element $C_i$ in the union. Its image under $\nu_{\infty}$ will be an irreducible component of $\pi^{-1}(\sing(X))$. This is so because there is an open subset in $C_i$ the image of which contains no arc which is the image of an arc in another component. 
\end{pf}
Next, we want to prove a similar claim for the $Y$-essential divisors over $X$. Let $\Xi$ be the disjoint union of the sets of $Y_i$-essential divisors for all $X_i$.
\begin{Prop}
The set of essential divisors over $X$ is in a one-to-one correspondence with $\Xi$.
\end{Prop}
\begin{pf}
Let $f': X' \rightarrow X$ be any $Y$-resolution of $X$ (that is, $f'$ is a divisorial resolution). Then by the universal property of the normalization map, $f'$ factors through $\nu$. If for each $i$, $f'_i: X'_i \rightarrow X_i$ is a $Y_i$-resolution of $X_i$, then $\coprod_i f'_i$ will be a resolution of $\coprod_i X_i$. So it defines a birational morphism $\psi: X' \rightarrow \coprod X'_i$. Take an essential divisor $D \subset X'$. The closure of $\psi(D)$ will then give an essential divisor over $X_i$ for some $i$. This is so because $\nu \circ (\coprod_i f'_i)$ is a $Y$-resolution for $X$, defining an element of $\Xi$. Thus, one defines an injective map on the set of essential divisors of $X$. This map is also surjective because the restriction of $f'$ to the pre-image of $(X_i) \subset X$ is a $Y_i$-resolution for $X_i$.
\end{pf}
By Prop.5.1 and Prop.5.2 we see that to prove the bijection of ${\mathcal{N}}_X$ it is enough show that the set $\Omega$ is in bijection with the set $\Xi$. But this follows from Thm.4.1 applied to each $X_i$, and the remark at the beginning of this section. This gives a positive answer for Nash problem in the case of STVs:
\begin{thm}
For an equidimensional STV $X$ there is a bijection between the set of the irreducible components of $\pi^{-1}(\sing(X))$ and the set of essential divisors over $X$.
\end{thm}

\bibliographystyle{alpha}
\bibliography{nash}
  
\end{document}